\documentclass[12pt]{amsart}
\usepackage{amssymb,amsmath}
\theoremstyle{plain}
\newtheorem{thm}{Theorem}[section]
\newtheorem{cor}[thm]{Corollary}
\newtheorem{prop}[thm]{Proposition}
\newtheorem{lem}[thm]{Lemma}

\theoremstyle{definition}
\newtheorem{rem}[thm]{Remark}
\newtheorem{defn}[thm]{Definition}

\newtheorem{eg}[thm]{Example}

\newcommand{\Prf}{\noindent\textbf{Proof.\ }}
\renewcommand{\qed}{\hfill \vrule height5pt width5pt depth1pt \vspace{+2.00ex}}
\newcommand{\upqed}{\vspace{-2.5\baselineskip}\newline\hbox{}\qed}

\newcommand{\bsl}{\setminus}

\newcommand{\ca}{\mathrm{C}^*}
\newcommand{\diag}{\operatorname{diag}}

\newcommand{\ep}{\varepsilon}
\renewcommand{\phi}{\varphi}
\newcommand{\lip}{\langle}
\newcommand{\rip}{\rangle}
\newcommand{\ip}[1]{\lip #1 \rip}
\newcommand{\mt}{\emptyset}
\newcommand{\ol}{\overline}
\newcommand{\one}{{\boldsymbol{1}}}
\newenvironment{sbmatrix}{\left[\begin{smallmatrix}}{\end{smallmatrix}\right]}

\newcommand{\sgn}{\operatorname{sgn}}

\newcommand{\upchi}{{\raise.35ex\hbox{$\chi$}}}
\renewcommand{\l}{\ell}
\newcommand{\Or}{\text{ or }}
\newcommand{\qand}{\quad\text{and}\quad}
\newcommand{\qfor}{\quad\text{for}\quad}
\newcommand{\qforal}{\quad\text{for all}\quad}
\newcommand{\qif}{\quad\text{if}\quad}

\newcommand{\foral}{\text{ for all }}
\newcommand{\Dim}{\operatorname{dim}}

\newcommand{\spn}{\operatorname{span}}
\newcommand{\Tr}{\operatorname{Tr}}
\newcommand{\bC}{{\mathbb{C}}}
\newcommand{\bN}{{\mathbb{N}}}

\newcommand{\bT}{{\mathbb{T}}}

\newcommand{\bZ}{{\mathbb{Z}}}
  \newcommand{\B}{{\mathcal{B}}}

  \newcommand{\G}{{\mathcal{G}}}
\renewcommand{\H}{{\mathcal{H}}}

\renewcommand{\P}{{\mathcal{P}}}

\renewcommand{\S}{{\mathcal{S}}}
  \newcommand{\T}{{\mathcal{T}}}
  \newcommand{\U}{{\mathcal{U}}}

  \newcommand{\X}{{\mathcal{X}}}
  
\newcommand{\fM}{{\mathfrak{M}}}

\newcommand{\fs}{{\mathfrak{s}}}
\newcommand{\fS}{{\mathfrak{S}}}

\begin{document}
\title{Norms of Schur Multipliers}
\author[K.R.Davidson]{Kenneth R. Davidson}
\address{Pure Math.\ Dept.\\U. Waterloo\\Waterloo, ON\;
N2L--3G1\\CANADA}
\email{krdavidson@math.uwaterloo.ca}
\author[A.P. Donsig]{Allan P. Donsig}
\address{Math.\ Dept.\\University of Nebraska\\
Lincoln, NE 68588\\USA}
\email{adonsig@math.unl.edu}
\date{}
\begin{abstract}
A subset $\P$ of $\bN^2$ is called Schur bounded if every infinite
matrix with bounded entries which is zero off of $\P$ yields a
bounded Schur multiplier on $\B(\H)$.
Such sets are characterized as being the union of a subset with
at most $k$ entries in each row with another that has at most $k$
entries in each column, for some finite $k$.
If $k$ is optimal, there is a Schur multiplier supported on the
pattern with norm $O(\sqrt k)$, which is sharp up to a constant.

The same techniques give a new, more elementary proof of
results of Varopoulos and Pisier on Schur multipliers with
given matrix entries of random sign.

We consider the Schur multipliers for certain matrices which have
a large symmetry group.  In these examples, we are able to compute
the Schur multiplier norm exactly.  This is carried out in detail for a
few examples including the Kneser graphs.
\end{abstract}
\maketitle

Schur multiplication is just the entrywise multiplication of matrices
or operators in a fixed basis.
These maps arise naturally as the (weak-$*$ continuous) bimodule maps
for the algebra of diagonal matrices (operators).
They are well-behaved completely bounded maps
that play a useful role in the theory of operator algebras.

As in the case of operators themselves, the actual calculation
of the norm of any specific Schur multiplier is a delicate task;
and is often impossible.  This has made it difficult to attack 
certain natural, even seemingly elementary, questions.

This study arose out of an effort to understand norms of Schur
multipliers supported on certain patterns of matrix entries.
The question of which patterns have the property that every
possible choice of bounded entries supported on the pattern yield
bounded Schur multipliers was raised by Nikolskaya and Farforovskaya 
in \cite{NF}.
We solve this problem completely.
The answer is surprisingly elegant.  The pattern must decompose
into two sets, one with a bound on the number of entries in each row,
and the other with a bound on the number of entries in each column.

There is a close relationship of these results with work of
Varopoulos \cite{Va} and Pisier \cite{Pi2}.  We had overlooked
this work and only discovered it late in our study.  Perhaps this
is just as well, as we may well have stopped had we realized
how close their results were to the ones we were seeking.
The upshot is that we also obtain a much more elementary
proof of the bulk of their results, though without the probabilistic
component.  Indeed our main tool in the decomposition is an
elementary, albeit powerful, combinatorial result known as the
min-cut-max-flow theorem.

In Section~\ref{S:HankToep}, we recover results of \cite{NF} on 
patterns of Hankel
and Toeplitz forms.  Actually the Toeplitz case is classical,
and we compare the bounds from our theorem with the
tighter bounds available from a deeper use of function theory.

Sections~\ref{S:symmetry} and \ref{S:graphs} deal with
exact computation of the Schur norm of certain matrices
that have lots of symmetry.  More precisely, let $G$ be a finite
group acting transitively on a set $X$.  
We obtain an explicit formula for the Schur multiplier norm of
matrices in the commutant of the action, i.e., matrices constant
on each orbit of $G$.
This uses a result of Mathias \cite{Ma}.
We carry this out for one nontrivial case---the adjacency 
matrix of the Kneser graph $K(2n+1,n)$, which has
$\binom{2n+1}{n}$ vertices indexed by $n$-element subsets
of $2n+1$, with edges between disjoint sets.

We would like to thank many people with whom we had helpful
conversations about aspects of this problem.  We thank Bill
Cunningham for showing us how to use the min-cut-max-flow theorem.
Thanks to Chris Godsil for pointing us to the literature on the
spectral analysis of Kneser and Johnson graphs.  We thank David
Handelman for sharing his notes on Johnson graphs with us.
We thank Stanislaw Szarek for pointing out a very useful paper 
of Fran\c coise Lust-Piquard \cite{LP};
and Nico Spronk for pointing out the paper by Pisier on
multipliers of nonamenable groups \cite{Pi2}.
Finally, we thank Vern Paulsen for some helpful comments.

\section{Background}\label{S:back}

If $A=[a_{ij}]_{i,j\in S}$ is a finite or infinite matrix,
the Schur (a.k.a.\ Hada\-mard) multiplier is the operator
$S_A$ on $\B(l^2(S))$ that acts
on an operator $T=[t_{ij}]$ by pointwise multiplication: 
$S_A(T) =[a_{ij}t_{ij}]$. 
To distinguish from the norm on bounded operators,
we will write $\|A\|_m$ for the norm of a Schur multiplier.
In general it is very difficult to compute the norm of a
Schur multiplier.  
Nevertheless, much is known in a theoretical
sense about the norm.  
In this section, we will quickly review some of the
most important results.

The following classical result owes most credit to Grothendieck.  
For a proof, see the books by Pisier \cite[Theorem 5.1]{Pi} and
Paulsen \cite[Theorem 8.7]{Pa}.

\begin{thm}\label{Schur_classic}
For $X$ an arbitrary set, let $S = [s_{ij}]$ be an $|X|\times|X|$ 
matrix with bounded entries considered as a Schur multiplier on $\B(l^2(X))$.
Then the following are equivalent:
\begin{enumerate}
\item $\|S\|_m \le 1$.
\item $\|S\|_{cb} \le 1$.
\item[($2'$)] There are contractions $V$ and $W$ from $l^2(X)$ to 
$l^2(X) \otimes l^2(Y)$ such that $S(A) = W^*(A \otimes I)V$.
\item There are unit vectors $x_i$ and $y_j$ in $l^2(Y)$
 so that $s_{ij} = x_i^*y_j$.
\item $\gamma_2(S) \le 1$ where
 $\gamma_2(S) = \inf_{S=AB} \|A\|_{2,\infty} \|B\|_{1,2}$.
\item There are $|X|\times|X|$ matrices $A=[a_{ij}]$ and
$B=[b_{ij}]$ with $a_{ii}=b_{ii}=1$ so that
$\begin{bmatrix}A&S\\S^*&B\end{bmatrix}$ is positive semidefinite.
\end{enumerate}
\end{thm}

Recall that the complete bound norm of $S$ is the norm of the
inflation of $S$ acting on operators with operator entries.
The most elegant proof of (1) implies (2) is due to Smith \cite{Sm}.
The converse is trivial.
The equivalence of (2) and ($2'$) is Wittstock's Theorem for representing
completely bounded maps.

The equivalence of (1), (3) and (4) is due to Grothendieck.
(3) follows from ($2'$) by taking 
$y_j = (E_{1j}\otimes I)Ve_j$ and $x_i = (E_{1i}\otimes I)We_i$.
Conversely, (3) implies ($2'$) by taking $Ve_j = e_j\otimes y_j$
and $We_i = e_i\otimes x_i$.
This condition was rediscovered by Haagerup, and
became well-known as his observation.  
So we shall refer to these as the Grothendieck--Haagerup
vectors for the Schur multiplier.

The $\gamma_2$ norm is the optimal factorization through Hilbert space
of $S$ considered as a map from $l^1$ to $l^\infty$.
The norm $\|A\|_{2,\infty}$ is the maximum of the 2-norm of the rows;
while $\|B\|_{1,2}$ is the maximum of the 2-norm of the columns.
Thus (3) implies (4) follows from $A = \sum_i e_ix_i^*$ and
$B = \sum_j y_j e_j^*$. And this implication is reversible.

The equivalence of (5) is due to Paulsen, Power and Smith \cite{PPS}.
This follows from (3) by taking $a_{ij} = x_i^*x_j$ and $b_{ij}=y_i^*y_j$.
Conversely, assume first that $X$ is finite.
Then the positive matrix $P$ decomposes as a sum of positive rank one matrices,
and thus have the form $[\bar{z}_iz_j]$ which can be seen to be a 
scalar version of (3).  Indeed it is completely positive. 
Hence the sum $S_P$ is also a completely positive  Schur multiplier.
Consequently $\|S_P\|_{cb} = \|S_P(I)\| = \max\{a_{ii}, b_{ii}\} = 1$.
So (2) holds. The case of general $X$ is a routine limit argument.

The $\gamma_2$ norm is equivalent to the norm in the Haagerup
tensor product $\l^\infty(X) \otimes_h \l^\infty(X)$, where we
identify an elementary tensor $a\otimes b$ with the matrix
$[a_ib_j]$. The Haagerup norm of a tensor $\tau$ is given by 
taking the infimum over all representations
$\tau = \sum_k a_k\otimes b_k$ of
\[
 \Big\| \sum_k a_ka_k^* \Big\|^{1/2} \, \Big\|\sum_k b_k^*b_k \Big\|^{1/2} .
\]
See \cite[Chapter 17]{Pa}.
Of course, since $\l^\infty$ is abelian, the order of the adjoints is
irrelevant.  
One can see the equivalence by taking a factorization
$S = AB$  from (4).  
Consider $A$ as a matrix with columns $a_k \in \l^\infty(X)$
and $B$ as a matrix with rows $b_k \in \l^\infty(X)$.
Identify the product with the tensor $ \sum_k a_k\otimes b_k$.
The norm $\|\sum_k a_ka_k^*\|^{1/2}$ can be seen to be $\|A\|_{2,\infty}$
and the norm of $\|\sum_k b_k^*b_k\|^{1/2}$ to be $\|B\|_{1,2}$.

\smallbreak
Generally, it is difficult to compute the norm of a Schur multiplier.
The exception occurs when the matrix $S$ is positive definite.
Then it is a classical fact that $S$ is a completely positive map.
Consequently, $\|S\|_{cb} = \|S(I)\| = \sup_{i\in X} s_{ii}$.

Grothendieck proved another remarkable result about Schur multipliers.
Recall that the projective tensor product $\l^\infty(X) \hat{\otimes}
\l^\infty(X)$ norms a tensor $\tau$ as the infimum over representations
$\tau = \sum_k a_k\otimes b_k$ of the quantity $\sum_k \|a_k\|\,\|b_k\|$.
It is a surprising fact that this norm and the $\gamma_2$ or Haagerup norm
are equivalent.  We will need this connection to understand the relevance
of work of Varopoulos.  For the moment, we state this result in a way that
makes a stronger connection to Schur multipliers.
An elementary tensor $a\otimes b$ yields a rank one matrix $[a_ib_j]$.
Thus Grothendieck's result is equivalent to:

\begin{thm}[Grothendieck] \label{Groth}
The convex hull of the rank one Schur multipliers of norm one
contains the ball of all Schur multipliers of norm at most $K_G^{-1}$,
where $K_G$ is a universal constant.
\end{thm}

In terms of the projective tensor product norm for a tensor $\tau$
and the corresponding Schur multiplier $S_\tau$, this result says that
\[
 K_G^{-1} \|\tau\|_{\l^\infty(S) \hat{\otimes} \l^\infty(S)} \le \|S_\tau\|_m
 \le \|\tau\|_{\l^\infty(S) \hat{\otimes} \l^\infty(S)}
\]
The constant $K_G$ is not known exactly.  In the complex case
Haagerup \cite{Hg2} showed that $1.338 < K_G < 1.405$;
and in the real case Krivine \cite{Kr} obtained the range $[1.676,1.783]$
and conjectured the correct answer to be $\frac\pi{2\log(1+\sqrt2)}$.

We turn to the results of Varopoulos \cite{Va} and Pisier \cite{Pi2}
which relate to our work.
The paper of Varopoulos is famous for showing that three commuting
contractions need not satisfy the von Neumann inequality.
Proofs of this, including the one in the appendix of Varopoulos's paper,
are generally constructive.
But the argument in the main part of his paper instead establishes
a result about $\l^\infty(X) \hat{\otimes} \l^\infty(X)$.
He does not establish precise information about constants.
This result was extended and sharpened by Pisier,
who casts it in  the language of Schur multipliers,
to deal with multipliers and lacunary sets on nonamenable groups.

Consider $\{\pm1\}^{X\times X}$ to be the space of functions
from $X\times X$ to $\{1,-1\}$ with the product measure
$\mu$ obtained from $p(1)=p(-1)=.5$.

\begin{thm}[Varopoulos--Pisier] \label{VP}
Let $S = [s_{ij}]$.  The following are equivalent.
\begin{enumerate}
\item For all $\ep\in \{\pm1\}^{X\times X}$,
 $\|[\ep_{ij}s_{ij}]\|_m < \infty$.
\item For almost all $\ep\in \{\pm1\}^{X\times X}$,
 $\|[\ep_{ij}s_{ij}]\|_m < \infty$.
\item $S = A+B$ and there is a constant $M$ so that 
\[ \sup_i \sum_j |a_{ij}|^2 \le M^2 \qand \sup_j \sum_i |b_{ij}|^2 \le M^2 .\]
\item There is a constant $M$ so that for every pair of finite subsets
$R,C \subset X$, $\sum_{i\in R,j\in C} |s_{ij}|^2 \le M^2 \max\{|R|, |C| \}$.
\end{enumerate}
\end{thm}

Pisier shows that if the average Schur multiplier norm
\[ \int \|[\ep_{ij}s_{ij}]\|_m \,d\mu(\ep) \le 1 ,\]
then one can take $M=1$ in (3).
Our results are not quite so sharp, as we require a 
constant (Lemma~\ref{bignorm}) of approximately 1/4.
The constant $M$ in the two conditions (3) and (4)
are not the same.  The correct relationship replaces
$\max\{|R|, |C| \}$ by $|R|+|C|$ (see Lemma~\ref{mincutmaxflow});
but they are related within a constant.
If $M$ is the bound in (3), it is not difficult to obtain a bound
of $2M$ for (1) (see Corollary~\ref{rowboundcor}).
Thus one obtains that the average Schur norm is within a factor
of 2 of the maximum.

\section{Schur Bounded Patterns} \label{S:patterns}

A \textit{pattern} $\P$ is a subset of $\bN \times \bN$.  
An infinite matrix $A = [a_{ij}]$ is \textit{supported on $\P$}
if $\{ (i,j) : a_{ij} \ne 0 \}$ is contained in $\P$.
We let $\S(\P)$ denote the set of Schur multipliers supported on
$\P$ with matrix entries $|s_{ij}|\le 1$.

More generally, we will also consider Schur multipliers dominated
by a given infinite matrix $A=[a_{ij}]$ with nonnegative entries.
Let $\S(A)$ denote the set of all Schur multipliers with
matrix entries $|s_{ij}|\le a_{ij}$.

\begin{defn}\label{d:schurbound}
Say that a pattern $\P \subset\bN \times \bN$ is \textit{Schur
bounded} if every $X \in \S(\P)$ yields bounded Schur
multiplier $S_X$.  
The Schur bound of $\P$ is defined as
$\fs(\P) := \sup_{X \in \S(\P)} \|X\|_m .$
Similarly, for a matrix $A$ with nonnegative entries, define
$\fs(\S(A)) = \sup_{X \in \S(A)} \|X\|_m$; and say that $\S(A)$
is Schur bounded if this value is finite.
\end{defn}

It is easy to see that if $\S(\P)$ is Schur bounded, 
then $\fs(\P)$ is finite.
Note that if $A_\P$ is the matrix with 1s on the entries of $\P$
and 0s elsewhere, then $\S(A_\P) = \S(\P)$.
We will maintain a distinction because we will require
\textit{integral decompositions} when working with a pattern $\P$.

Certain patterns are easily seen to be Schur bounded and this
is the key to our result.
The following two definitions of row bounded for patterns and matrices 
are not parallel, as the row bound of $A_\P$ is actually 
the square root of the row bound of $\P$.
Each definition seems natural for its context, so we content ourselves
with this warning.

\begin{defn}
A pattern is \textit{row bounded by $k$} if there are at most $k$ 
entries in each row; and \textit{row finite} if it is row bounded
by $k$ for some $k\in\bN$.  Similarly we define \textit{column
bounded by $k$} and \textit{column finite}. 

A nonnegative matrix $A=[a_{ij}]$ is \textit{row bounded by $L$}
if the rows of $A$ are bounded by $L$ in the $l^2$-norm: 
$\sup_{i\ge1} \sum_{j\ge1} |a_{ij}|^2 \le L^2 < \infty$.
Similarly we define \textit{column bounded by $L$}.
\end{defn}

The main result of this section is:

\begin{thm}\label{Schurbound}
For a pattern $\P$, the following are equivalent:
\begin{enumerate}
\item $\P$ is Schur bounded.
\item $\P$ is the union of a row finite set and a column finite set.
\vspace{.5ex}
\item $\displaystyle\sup_{R,C \text{ finite}}\frac{|\P \cap (R\times C)|}
{|R|+|C|} < \infty$.
\end{enumerate}
Moreover, the optimal bound $m$ on the size of the row and column finite sets
in $(2)$ coincides with the least integer dominating the supremum in $(3)$; and
the Schur bound satisfies
\[ \sqrt m /4 \le \fs(\P) \le 2\sqrt m .\]
\end{thm}

This theorem has a direct parallel for nonnegative matrices.

\begin{thm}\label{matrix_version}
For a nonnegative infinite matrix $A = [a_{ij}]$, the following are equivalent:
\begin{enumerate}
\item $\S(A)$ is Schur bounded.
\item $A=B+C$ where $B$ is row bounded and $C$ is column bounded.
\item $\displaystyle\sup_{R,C \text{ finite}} 
 \frac{\sum_{i\in R, j \in C}a_{ij}^2}{|R|+|C|} < \infty$.
\end{enumerate}
Moreover, the optimal bound $M$ on the row and column bounds in $(2)$
coincides with the square root of the supremum $M^2$ in $(3)$; 
and the Schur bound satisfies
\[ M /4 \le \fs(\P) \le 2M .\]
\end{thm}

\begin{lem}\label{rowbound}
If $\P$ is row $($or column$)$ bounded by $n$, 
then $\fs(\P) \le \sqrt n$.

Likewise if $A$ is row $($or column$)$ bounded by $L$, then
$\fs(\S(A)) \le L$.
\end{lem}

\Prf 
The pattern case follows from the row bounded case for the
nonnegative matrix $A = A_\P$ with $L = \sqrt n$.
Suppose that $\S(A)$ is row bounded by $L$.
Consider any $S \in \S(A)$.
Then $\sup_{i\ge1} \sum_{j\ge1} |s_{ij}|^2 \le L^2$.
Define vectors $x_i = \sum_{j\ge1} s_{i,j} e_j$ for $i \ge1$.
Then $\sup_{i\ge1} \|x_i\| \le L$;
and $\ip{x_i,e_j} = s_{ij}$.
So by the Grothendieck--Haagerup condition,
\[
 \| S \|_m \le \sup_{i,j} \|x_i\| \, \|e_j\| \le L.
\]
Thus $\fs(\S(A)) \le L$.
\qed

\begin{cor}\label{rowboundcor}
If $\P$ is the union of a set row bounded by $n$
and a set column bounded by $m$, then
$\P$ is Schur bounded with bound $\sqrt n + \sqrt m$.

Likewise, if $A = B + C$ such that $B$ is row bounded by
$L$ and $C$ is column bounded by $M$, then $\fs(\S(A)) \le L+M$.
\end{cor}

We require a combinatorial characterization of sets which are the
union of an $n$-row bounded set and an $m$-column bounded set.
This will be a consequence of the min-cut-max-flow theorem
(see \cite{CCPS}, for example).  
This is an elementary result in combinatorial
optimization that has many surprising consequences.  
For example, it has been used by Richard Haydon to give a short proof 
of the reflexivity of commutative subspace lattices~\cite{Ha}.
It should be more widely known.

\begin{lem}\label{mincutmaxflow}
A pattern $\P$ is the union of a set $\P_r$ row bounded by $m$
and a set $\P_c$ column bounded by $n$ if and only if for every 
pair of finite subsets $R,C \subset \bN$, 
\[ | \P \cap R \times C | \le m |R| + n |C| . \]

Similarly, a matrix $A = [a_{ij}]$ with nonnegative entries
decomposes as a sum $A = A_r+A_c$ where $A_r$ is row bounded by $M^{1/2}$
and $A_c$ is column bounded by $N^{1/2}$ if and only if for every 
pair of finite subsets $R,C \subset \bN$, 
\[ \sum_{i\in R}\sum_{j\in C} a_{ij}^2 \le M |R| + N |C| . \]
\end{lem}

\Prf The two proofs are essentially identical.
However the decomposition of $\P$ must be into two disjoint subsets.
This means that the decomposition $A_\P = A_{\P_1} + A_{\P_2}$
is a split into $0,1$ matrices.
We will work with $A$, but will explain the differences in the
pattern version when it arises.

The condition is clearly necessary.

For the converse, we first show that it suffices to solve
the finite version of the problem.
For $p \in \bN$, let $A_p$ be the restriction of $A$ to the
first $p$ rows and columns.
Suppose that we can decompose $A_p = A_{r,p} + A_{c,p}$
where $A_{r,p}$ is row bounded by $M^{1/2}$
and $A_{c,p}$ column bounded by $N^{1/2}$ for each $p \in \bN$.
Fix $k$ so that $A_k \ne 0$.  
For each $p\ge k$, the set of such decompositions
for $A_p$ is a compact subset of $\fM_p \times \fM_p$.
In the pattern case, we consider only $0,1$ decompositions.
The restriction to the $k\times k$ corner
is also a compact set, say $\X_{k,p}$.
Observe that this is a decreasing sequence of nonempty compact sets.
Thus $\cap_{p\ge k} \X_{k,p} = \X_k$ is nonempty.
Therefore there is a consistent choice of a decomposition 
$A = A_r+A_c$ so that the restriction to each $k\times k$ 
corner lies in $\X_k$ for each $k\ge1$.
In the pattern case, the entries are all zeros and ones.

So now we may assume that $A=[a_{ij}]$ is a matrix supported on $R_0 \times C_0$,
where $R_0$ and $C_0$ are finite.
We may also suppose that the $l^2$-norm of each row is greater than $M^{1/2}$
and the $l^2$-norm of each column is greater than $N^{1/2}$.
For otherwise, we assign all of those entries in the row to $A_r$ 
(or all entries in the column to $A_c$) and delete the row (column).
Solving the reduced problem will suffice.
If after repeated use of this procedure, the matrix is empty, we are done.
Otherwise, we reach a reduced situation in which the $l^2$-norm of each      
row is greater than $M^{1/2}$ and the $l^2$-norm of each column is greater  
than $N^{1/2}$.

Define a graph $\G$ with vertices $\alpha$, $r_i$ for $i \in R_0$,
$c_j$ for $j \in C_0$, and $\omega$.  
Put edges from each $r_i \in R_0$ to each $c_j \in C_0$,
from $\alpha$ to $r_i$, $i \in R_0$,
and from $c_j$ to $\omega$, $j \in C_0$.
Consider a network flow on the graph in which the edge from 
$r_i$ to $c_j$ may carry $a_{ij}$ units;
edges leading out of $\alpha$ can carry up to $M$ units;
and the edge from $c_j$ to $\omega$ can carry $v_j-N$ units,
where $v_j = \sum_{i\in R_0} a_{ij}^2$. 
In the pattern case, these constraints are integers.

The min-cut-max-flow theorem states that the maximal possible flow
from $\alpha$ to $\omega$ across this network equals the minimum flow
across any cut that separates $\alpha$ from $\omega$.
Moreover, when the data is integral, the maximal flow comes from an integral
solution. 
A \textit{cut} $\X$ is just a partition of the vertices into two
disjoint sets
$\{\alpha\} \cup R_1 \cup C_1$ and $\{\omega\} \cup R_2 \cup C_2$.
The flow across the cut is the total of allowable flows on each edge
between the two sets.

The flow across the cut $\X$ is
{\allowdisplaybreaks
\begin{align*}
 f(\X) &= \sum_{i\in R_1}\sum_{j\in C_2} a_{ij}^2 +M |R_2| 
              + \sum_{j\in C_1} (v_j - N) \\
 &=  \sum_{i\in R_1}\sum_{j\in C_2} a_{ij}^2 + M |R_2| - N |C_1|
      + \sum_{i\in R_0}\sum_{j\in C_1} a_{ij}^2 \\
 &= \sum_{i\in R_0}\sum_{j\in C_0} a_{ij}^2 - 
    \sum_{i\in R_2}\sum_{j\in C_2} a_{ij}^2 + M |R_2| + N |C_2| - N |C_0| \\
  &\ge \sum_{i\in R_0}\sum_{j\in C_0} a_{ij}^2 - N |C_0|
\end{align*}
}
The last inequality uses the hypothesis on $A$ with $R=R_2$ and $C=C_2$.
On the other hand, the cut separating $\omega$ from the rest has flow
exactly 
\[
 \sum_{j\in C_0} (v_j - N) = \sum_{i\in R_0}\sum_{j\in C_0} a_{ij}^2 - N |C_0| .
\]

Therefore there is a network flow that achieves this maximum.
In the pattern case, the solution is integral.
Necessarily this will involve a flow of exactly $v_j - N$ from
each $j \in C_0$ to $\omega$. 
Let $b_{ij}$ be the optimal flow from $r_i$ to $c_j$.
So $0 \le b_{ij} \le a_{ij}$.
The flow out of each $r_i$ equals the flow into $r_i$ from $\alpha$,
whence $\sum_{j\in C_0} b_{ij} \le M$.

Define the matrix $A_r = \big[ \sqrt{b_{ij}} \big]$ and
$A_c = \big[ \sqrt{a_{ij}-b_{ij}} \big]$.
In the pattern case, these entries are 0 or 1.
Then the rows of $A_r$ are bounded by $M^{1/2}$.
The $j$th column of $A_c$ has norm squared equal to
\[
 \sum_{i\in R_0} a_{ij} - b_{ij} = v_j - (v_j-N) = N .
\]
This is the desired decomposition and it is integral for patterns.
\qed

To construct large norm Schur multipliers on certain patterns, 
we will make use of the following remarkable result by Fran\c coise
Lust-Piquard \cite[Theorem~2]{LP}. 
While the method of proof is unexpected, it is both short and elementary.

\begin{thm}[Lust-Piquard]\label{LPthm}
Given any $($finite or infinite$)$ nonnegative matrix $X = [x_{ij}]$ satisfying
\[
 \max_i \sum_j x_{ij}^2 \le 1 \qand \max_j \sum_i x_{ij}^2 \le 1
 \qforal i,j,
\]
there is an operator $Y = [y_{ij}]$ so that 
\[
 \|Y\| \le \sqrt6 \qand |y_{ij}| \ge x_{ij} \foral i,j.
\]
\end{thm}

The constant of $\sqrt{6}$ is optimal, as shown in an addendum to~\cite{LP}.

\begin{lem}\label{bignorm}
Let $A =[a_{ij}]$ be a nonnegative $m \times m$ matrix
such that $\sum_{i=1}^m \sum_{j=1}^m a_{ij}^2 = m\alpha$.\vspace{.5ex}
Then there is a Schur multiplier $S \in \S(A)$ 
such that $\|S\|_m \ge \frac12 \sqrt{\frac \alpha 3}$.
\end{lem}

\Prf We may assume that there are no nonzero rows or columns.
Let 
\[
 r_i = \sum_{j=1}^m a_{ij}^2 \qand c_j = \sum_{i=1}^m a_{ij}^2 .
\]
Define 
\[
 x_{ij} = \frac{a_{ij}}{\sqrt{r_i+c_j}} .
\]
Let $X = [x_{ij}]$.
The row norms of $X$ satisfy
\[
 \sum_{j=1}^m x_{ij}^2 \le
 \sum_{j=1}^m \frac{a_{ij}^2}{r_i} = 1 ;
\]
and similarly the column norms are bounded by 1.

By Theorem~\ref{LPthm}, there is a matrix $Y$ such that
\[
 \|Y\| \le \sqrt6 \qand |y_{ij}| \ge x_{ij} \qforal i,j. 
\]
Define $s_{ij} = a_{ij}x_{ij}/y_{ij}$ (where $0/0 :=0$).
Then $S = [s_{ij}]$ belongs to $\S(A)$.
Observe that 
\[
 S(Y) = Z := [a_{ij}x_{ij}]
  = \begin{bmatrix} \frac{a_{ij}^2}{\sqrt{r_i+c_j}}\end{bmatrix} .
\]  
Hence $\|S\|_m \ge \|Z\|/K$.

Define vectors $u = (u_i)$ and $v=(v_j)$ by 
\[
 u_i = \Big(\frac{r_i}{m\alpha}\Big)^{1/2} \qand 
 v_j = \Big(\frac{c_j}{m\alpha}\Big)^{1/2}.
\]
Then $\|u\|_2^2 = \displaystyle\frac1{m\alpha} \sum_{i=1}^m r_i = 1$
and similarly $\|v\|_2=1$.
Compute 
\[
  \|Z\| \ge u^* Z v
  = \frac1{m\alpha} \sum_{i=1}^m \sum_{j=1}^m a_{ij}^2
     \sqrt{\frac{r_ic_j}{r_i+c_j}} .
\]
Observe that $\displaystyle\sqrt{\frac{r_ic_j}{r_i+c_j}} =
 \Big(\frac1{r_i} + \frac1{c_j} \Big)^{-1/2}$.
Also
\begin{align*}
 \sum_{i=1}^m \sum_{j=1}^m a_{ij}^2 \Big(\frac1{r_i} + \frac1{c_j} \Big) &=
 \sum_{i=1}^m \sum_{j=1}^m \frac{a_{ij}^2}{r_i}
  + \sum_{j=1}^m \sum_{i=1}^m \frac{a_{ij}^2}{c_j} \\
 &= \sum_{i=1}^m 1 + \sum_{j=1}^m 1 = 2m.
\end{align*}

A routine Lagrange multiplier argument shows that if $\alpha_k \ge 0$ are
constants, $t_k > 0$ are variables, and 
$\sum_{k=1}^{m^2} \alpha_k t_k = 2m$, then
$\sum_{k=1}^{m^2} \alpha_k t_k^{-1/2}$ is minimized
when all $t_k$ are equal.
Hence if $\sum_{k=1}^{m^2} \alpha_k = m\alpha$,
\[
 \sum_{k=1}^{m^2} \alpha_k t_k^{-1/2}
  \ge m\alpha \Big(\frac{2m}{m\alpha}\Big)^{-1/2}
  = m\alpha \sqrt{\frac \alpha 2}
\]
Applying this to the numbers $\frac1{r_i} + \frac1{c_j} $ yields
\[
 \|Z\| \ge \frac1{m\alpha} \sum_{i=1}^m \sum_{j=1}^m a_{ij}^2
 \Big(\frac1{r_i} + \frac1{c_j} \Big)^{-1 /2} \ge\sqrt{\frac \alpha 2}.
\]
Thus $\|S\|_m \ge \frac{\sqrt{\alpha}}{\sqrt6 \sqrt2} = \frac12 \sqrt{\frac{\alpha}3}$.
\qed

\bigskip\pagebreak[3]
\noindent\textbf{Proof of Theorem~\ref{Schurbound} and
Theorem~\ref{matrix_version}.}

Statements (2) and (3) are equivalent by Lemma~\ref{mincutmaxflow},
taking $m=n$ and $M=N$.

Assuming (2), Corollary~\ref{rowboundcor} shows that $\P$ or $\S(A)$ is Schur
bounded by $2\sqrt m$ or $2M$. 
Assuming (3) in the pattern case, the supremum exceeds $m-1$; so
Lemma~\ref{bignorm} shows that 
\[
 \fs(\P) \ge \frac{\sqrt{m-1}}{2\sqrt3} \ge \frac{\sqrt m } 4
\]
for $m \ge 4$. 
For $m \le 16$, $\sqrt m /4 \le 1$; and $1$ is also a lower bound for any pattern. 
For the matrix case, we use the exact supremum in Lemma~\ref{bignorm}, so we
obtain a lower bound of $M/4$.

Conversely, if the supremum in (3) is infinite, the same argument shows that
the Schur bound is infinite.
In fact it is easy to see that this implies that $\S(\P)$ or $\S(A)$ contains
unbounded Schur multipliers.  
It is not difficult to produce disjoint finite rectangles $R_n \times C_n$ on
which the ratio in (3) exceeds $n^2$. 
So by Lemma~\ref{bignorm}, we construct a Schur multiplier $S_n$ in 
$\S(\P)$ or $\S(A)$ supported on $R_n \times C_n$ with Schur norm at least 
$n/4$.
Take $S$ to be defined on each rectangle as $S_n$ and zero elsewhere.
Then $S$ is an unbounded Schur multiplier in this class.
\qed

\begin{rem}
One might suspect, from the $\sqrt n$ arising in Lemma~\ref{rowbound}, that
if two matrices are supported on pairwise disjoint patterns, there might
be an $L^2$ estimate on the Schur norm of the sum.
This is not the case, as the following example shows.

Let $\one = (1,1,1,1)^t \in \bC^4$ and $A = \one\one^* - I$.
If $U=\diag(1,i,-1,-i)$, then the diagonal expectation is 
\[
 \Delta(X) = S_I(X) = \frac14 \sum_{k=0}^3 U^k X U^{*k} .
\]
We use a device due to Bhatia--Choi--Davis \cite{BCD}.  Observe that
\begin{align*}
 S_{A+tI}(X) &= X + (t-1) \Delta(X) \\
             &= (1 + \frac{t-1}4 )X + \frac{t-1}4 \sum_{k=1}^3 U^kXU^{*k} .
\end{align*}
Therefore
\begin{align*}
  \|S_{A+tI}\|_m &\le \Big| 1 + \frac{t-1}4 \Big|  + \frac{3|t-1|}4  
  \\&= \begin{cases}
      |t|       &\qif t \ge 1 \Or t \le -3\\
      \frac12 |3-t|   &\qif -3 \le t \le 1
\end{cases} .
\end{align*}
On the other hand, $S_{A+tI}(I) = tI$; so $\| S_{A+tI} \|_m \ge |t|$.
Observe that $\frac14 \one\one^*$ is a projection.
Hence $A+tI = \one\one^* + (t-1)I$ has spectrum $\{t-1,t+3\}$;
and thus 
\[ \|A+tI\| = \max\{|t-1|, |t+3|\} . \]
So $\|A-I\| = 2$.
If $-3 \le t \le 1$, then $S_{A+tI}(A-I) = A-tI$ has norm $|3-t|$
and so $\|S_{A+tI}\|_m \ge |3-t|/2$.

In particular, $\|S_A\|_m = \frac32$ and $\|S_I\|_m=1$, but
\[ \|S_{A-I}\|_m = 2 > \big( \|S_A\|_m^2 + \|S_I\|_m^2 \big)^{1/2} .\]
\end{rem}

\begin{rem}
In \cite{BGN}, Bennett, Goodman and Newman show that if
$A$ is an $n \times n$ matrix with entries taking the values $\pm1$
with probability .5, then on average the norm of $A$ is bounded
by $K\sqrt n$, where $K$ is a universal constant.
This is best possible as each row and column has norm $\sqrt n$.
This minimum can be achieved in certain cases, for example
by tensoring copies of 
$\begin{sbmatrix}\phantom{-}1&1\\-1&1\end{sbmatrix}$ together.  
The maximum norm occurs for the matrix $\one\one^*$
for which all entries are 1, in which case the norm is $n$.
So we see that the average norm is within
a constant of the minimum.

This can be used to show that, on average, the Schur norm
$\|A\|_m$ is near the maximum $\sqrt n$.  
Indeed, $S_A(A) = \one\one^*$.
So 
\[ \|A\|_m \ge \frac n {\|A\|} \ge K^{-1} \sqrt n \]
on average.
\end{rem}

\section{Hankel and Toeplitz Patterns} \label{S:HankToep}

A \textit{Hankel pattern} is a set of the form
\[
 \H(S) = \{(i,j) : i,j \in \bN,\ i+j \in S \} \qfor S \subset \bN .
\]
A \textit{Toeplitz pattern} is a set of the form
\[
 \T(S) = \{ (i,j) : i,j \in \bN_0,\ i-j \in S \} \qfor S \subset \bZ. 
\]  
Recall that a set $S = \{ s_1 < s_2 < \dots \}$ is lacunary if there
is a constant $q>1$ so that $s_{i+1}/s_i > q$ for all $i \ge1$.

Nikolskaya and Farforovskaya show that a Hankel pattern is Schur bounded 
if and only if it is a finite union of lacunary sets~\cite[Theorem~3.8]{NF},
by considering Fej\'er kernels and Toeplitz extensions.
We give an elementary proof based on Theorem~\ref{Schurbound}.

\begin{prop}\label{Hankel}
Consider a Hankel pattern $\H(S)$ of a set $S \subset \bN$.
Then the following are equivalent:
\begin{enumerate}
\item $\H(S)$ is Schur bounded.
\item $\H(S)$ is the union of a row finite and a column finite set.
\item $\sup_{k\ge0} |S \cap (2^{k-1},2^k]| < \infty$.
\item $S$ is the union of finitely many lacunary sets.  
\end{enumerate}
\end{prop}

\Prf By Theorem~\ref{Schurbound}, (1) and (2) are equivalent.

Let $a_k = |S \cap (2^{k-1},2^k]|$ for $k \ge 0$.
If (3) holds, $\max_{k\ge0} a_k = L < \infty$.
So $S$ splits into $2L$ subsets with at most one element in every second
interval $(2^{k-1},2^k]$; which are therefore lacunary with ratio at least 2.
Conversely, suppose that $S$ is the union of finitely many lacunary sets.
A lacunary set with ratio $q$ may be split into $d$ lacunary sets of ratio 2
provided that $q^d\ge2$.  So suppose that there are $L$ lacunary sets of ratio
2. Then each of these sets intersects $(2^{k-1},2^k]$ in at most one element.
Hence $\max_{k\ge0} a_k \le L < \infty$.
Thus (3) and (4) are equivalent.

Suppose that $S$ is the union of $L$ sets $S_i$ which are each lacunary with
constant $2$. 
Split each $\H(S_i)$ into the subsets $R_i$ on or below the diagonal and
$C_i$ above the diagonal.  
Observe that $R_i$ is row bounded by 1, and $C_i$ is column bounded by 1.
Hence (4) implies (2).

Consider the subset of $\H(S)$ in the first $2^k$ rows and columns
$R_k \times C_k$.
This square will contain at least $2^{k-1}a_k$ entries corresponding
to the backward  diagonals for $S \cap (2^{k-1},2^k]$, which all have
more than $2^{k-1}$ entries.
Thus 
\[
 \sup_{k \ge 0} \frac{|\H(S) \cap (R_k\times C_k)|}{|R_k|+|C_k|}
 \ge \sup_{k \ge 0} \frac{2^{k-1}a_k}{2^k+2^k} = \sup_{k \ge 0} \frac{a_k}4 .
\]
Hence if (3) fails, this supremum if infinite.
Thus $\H(S)$ is not the union of a row finite and a column finite set.
So (2) fails.
\qed

The situation for Toeplitz patterns is quite different.
It follows from classical results, as we explain below, and 
Nikolskaya and Farforovskaya outline a related proof~\cite[Remark~3.9]{NF}.
But first we show how it follows from our theorem.

\begin{prop}\label{Toeplitz}
The Toeplitz pattern $\T(S)$ of any infinite set $S$ is not Schur bounded.
Further,  
\[ \frac14 |S|^{1/2} \le \fs(\T(S)) \le |S|^{1/2} .\]
\end{prop}

\Prf Since $\T(S)$ is clearly row bounded by $|S|$, the upper bound
follows from Lemma~\ref{rowbound}.

Suppose that $S = \{ s_1 < s_2 < \dots < s_n \}$.
Consider the $m \times m$ square matrix with upper left hand corner
equal to $(s_1,0)$ if $s_1\ge0$ or $(0,-s_1)$ if $s_1<0$.
Then beginning with row $m-(s_n-s_1)$, there will be $n$ entries of $\T(S)$ in
each row.  Thus the total number of entries is at least $n(m-(s_n-s_1))$.
For $m$ sufficiently large, this exceeds $(n-1)m$.
Hence by Lemma~\ref{bignorm}, 
\[
 \fs(\T(S)) \ge \frac{\sqrt{n-1}}{2\sqrt3} \ge \frac{\sqrt n}{4}
\]
provided $n\ge4$. The trivial lower bound of 1 yields the lower
bound for $n<4$. 
\qed

To see how this is done classically, 
we recall the following \cite[Theorem~8.1]{Be}.
Here, $\T$ denotes the space of Toeplitz operators.

\begin{thm}[\textbf{Bennett}]
A Toeplitz matrix $A = [a_{i-j}]$ determines a bounded Schur multiplier
if and only if there is a finite complex Borel measure $\mu$ on
the unit circle $\bT$ so that $\hat{\mu}(n) = a_n$, $n \in \bZ$.  Moreover
\[
 \|A\|_m = \| S_A|_{\T} \| = \| \mu \| .
 \]
\end{thm}

We combine this with estimates obtained from the Khintchine inequalities.

\begin{thm}\label{Toeplitz-Khinchine}
Let $(a_k)_{k\in\bZ}$ be an $l^2$ sequence and let $A = [a_{i-j}]$.  
Then
\[ \frac1{\sqrt2} \|(a_k)\|_2 \le \fs(A) \le \|(a_k)\|_2 .\]
\end{thm}

\Prf
Suppose $S \in \S(A)$, that is, $S = [s_{ij}]$ with $|s_{ij}| \le a_{i-j}$.
Then each row of $S$ has norm bounded by $\|(a_k)\|_2$.
Hence by Lemma~\ref{rowbound}, $\|S\|_m \le \|(a_k)\|_2$.
So $\fs(A) \le \|(a_k)\|_2$.

Conversely, let $X := \{1,-1\}^\bZ$.
Put the measure $\mu$ on $X$ which is the product of
measures on $\{-1,1\}$ assigning measure $1/2$ to both $\pm1$.
For $\ep = (\ep_k)_{k \in \bZ}$ in $X$,
define $f_\ep(\theta) = \sum_{k\in\bZ} \ep_k a_k e^{ik\theta}$.
Then $f_\ep \in L^2(\bT) \subset L^1(\bT)$.
Hence $S_\ep := S_{T_{f_\ep}}$ defines a bounded Schur multiplier
with 
\[
 \|S_\ep\|_m = \|f_\ep\|_1 \le \|f_\ep\|_2 = \|(a_k)\|_2 .
\]
Then we make use of the Khintchine inequality \cite{Sz,Hg}:
\[
 \frac1{\sqrt2} \|(a_k)\|_2 \le 
 \int_X \|f_\ep\|_1 \,d\mu(\ep) \le \|(a_k)\|_2 .
\]
It follows that on average, most $f_\ep$ have $L^1$-norm
comparable to the $L^2$-norm.  
In particular, there is some choice of $\ep$ with 
$\|f_\ep\|_1 \ge \frac1{\sqrt2} \|(a_k)\|_2$.
Thus $\fs(A) \ge \|S_\ep\|_m \ge \frac1{\sqrt2} \|(a_k)\|_2 .$
\qed

\begin{rem}\label{Toeplitz_rem}
In the case of a finite Toeplitz pattern $\T(S)$,
say $S = \{s_1 < s_2 < \dots < s_n \}$,
$f_\ep = \sum_{k=1}^n \ep_k e^{is_k\theta}$.
We can use the Khintchine inequality for $L^\infty$:
\[
 \|(a_k)\|_2 \le \int_X \|f_\ep\|_\infty \,d\mu(\ep)
 \le \sqrt2 \|(a_k)\|_2 .
\]
Thus there will be choices of $\ep$ so that 
$\|f_\ep\|_\infty \le \sqrt{2n}$.
Then note that $S_{T_{f_\ep}}(T_{f_\ep}) = T_{f_\one}$,
where $f_\one = \sum_{k=1}^n e^{is_k\theta}$.
Clearly $\|f_\one\|_\infty = f_\one(0) = n$.
Thus $\|S_{T_{f_\ep}}|_{\T(S)}\| \ge \sqrt{n/2}$.
\end{rem}

\section{Patterns with a Symmetry Group} \label{S:symmetry}

Consider a finite group $G$ acting \textit{transitively} on a finite set $X$.
Think of this as a matrix representation on the Hilbert space $\H_X$
with orthonormal basis $\{e_x : x \in X \}$.
Let $\pi$ denote the representation of $G$ on $\H_X$ and
$\T$ the commutant of $\pi(G)$.
The purpose of this section is to compute the norm
of $S_T$ for $T \in \T$.

Decompose $X^2$ into G-orbits $X_i$ for $0 \le i \le n$, 
beginning with the diagonal $X_0 = \{(x,x) : x \in X\}$. 
Let $T_i \in \B(\H_X)$ denote the matrix with
$1$s on the entries of $X_i$ and $0$ elsewhere.
Then it is easy and well-known that $\T$ is
$\spn \{T_i : 0 \le i \le n \}$.
In particular, $\T$ is a C*-algebra.
Also observe that every element of $\T$ is constant on the 
main diagonal.  

Since $G$ acts transitively on $X$, $r_i := |\{y \in X : (x,y) \in X_i \}|$
is independent of the choice of $x\in X$.
Thus the vector $\one$ of all ones is a common eigenvector
for each $T_i$, and hence for all elements of $\T$,
corresponding to a one-dimensional reducing subspace
on which $G$ acts via the trivial representation.

First we establish an easy, general upper bound for $\|T\|_m$
where $T \in \T$.
As usual, $\Delta$ is the expectation onto the diagonal.

\begin{prop}\label{schur_bound}
For a matrix $T$, 
\[
 \|T\|_m \le \| \Delta(|T^*|) \|^{1/2} \, \| \Delta(|T|) \|^{1/2}
 = \bigl\|\, |T^*|\, \bigr\|_m^{1/2} \, \bigl\|\,|T|\,\bigr\|_m^{1/2}
.\]
\end{prop}

\Prf Use polar decomposition to factor $T = U |T|$.
Define vectors $x_i = |T|^{1/2} e_i$ and $y_j = |T|^{1/2}U^* e_j$.
Then
\[
 \ip{x_i,y_j} = \ip{|T|^{1/2} e_i, |T|^{1/2}U^* e_j} = \ip{Te_i,e_j} .
\]
This yields a Grothendieck--Haagerup form for $S_T$. Now 
\[
 \|x_i\|^2 = \ip{|T|^{1/2} e_i, |T|^{1/2} e_i} = \ip{|T|e_i,e_i} .
\]
Hence $\max_i \|x_i\| = \| \Delta(|T|) \|^{1/2}$.
Similarly, since $|T|^{1/2}U^* = U^* |T^*|^{1/2}$
\[
 \|y_j\|^2 = \ip{U^* |T^*|^{1/2} e_j, U^* |T^*|^{1/2} e_j} 
 = \ip{|T^*| e_j,e_j}.
\]
So $\max_j \|y_j\| = \| \Delta(|T^*|) \|^{1/2}$.
Therefore
\[
 \|T\|_m \le \max_{i,j} \|x_i\|\,\|y_j\| = 
 \| \Delta(|T^*|) \|^{1/2} \, \| \Delta(|T|) \|^{1/2} .
\]
Since $|T|$ and $|T^*|$ are positive, the Schur norm is just
the sup of the diagonal entries.
\qed

\begin{cor}\label{schur_bound_sa}
If $T=T^*$, then $\|T\|_m \le \| \Delta(|T|) \|$.
\end{cor}

\begin{rem}
In general this is a strict inequality.  
If $T = \begin{bmatrix}4&3\\3&1\end{bmatrix}$, then
$|T|= \begin{bmatrix}2\sqrt5&\sqrt5\\ \sqrt5&\sqrt5\end{bmatrix}$.
But $\|S_T\|_m= 4 < 2\sqrt5$.  Indeed, take $x_1=y_1=2e_1$
and $x_2=\frac32 e_1 + \frac{\sqrt5}2 e_2$ and
$y_2 = \frac32 e_1 - \frac{\sqrt5}2 e_2$.
\end{rem}

The main result of this section is:

\begin{thm}\label{schur_group}
Let $X$ be a finite set with a transitive action by a finite group $G$.
If $T$ belongs to $\T$, the commutant of the action of $G$,
then for any $x_0\in X$,
\[
 \| T \|_m = \| S_T|_\T \|
 = |X|^{-1} \Tr(|T|) = \ip{|T|e_{x_0},e_{x_0}} .
\]
\end{thm}

This result is a special case of a nice result of Mathias \cite{Ma}.
As far as we know, the application of Mathias' result to
the case of matrices invariant under group actions has not
been exploited.
As Mathias's argument is short and elegant, we include it.

\begin{thm}[Mathias]\label{mathias}
If $T$ is an $n\times n$ matrix with
$\Delta(|T^*|)$ and $\Delta(|T|)$ scalar, then
\[
 \| T \|_m =  \frac1n  \Tr(|T|) .
\]
\end{thm}

\Prf 
For an upper bound, Proposition~\ref{schur_bound} shows that 
\begin{align*}
  \| T\|_m &\le \| \Delta(|T^*|) \|^{1/2} \, \| \Delta(|T|) \|^{1/2} \\
 &=  \big(\tfrac1n\Tr(|T^*|) \big)^{1/2}  \big(\tfrac1n \Tr(|T|) \big)^{1/2} 
  = \tfrac1n \Tr(|T|),
\end{align*}
because $|T|$ and $|T^*|$ are constant on the main diagonal,
and $|T^*|$ is unitarily equivalent to $|T|$, and so has the same trace.

For the lower bound, use the polar decomposition $T = W |T|$.
Let $\ol{W}$ have matrix entries which are the complex conjugates
of the matrix entries of $W$.
Write $T = [t_{ij}]$ and $W = [w_{ij}]$ as $n \times n$ matrices
in the given basis.
Set $\one$ to be the vector with $n$ 1's.
Then
\begin{align*}
 \| T \|_m &\ge \|S_T(\ol{W}) \|  
  \ge \frac 1 n \ip{ S_T(\ol{W}) \one, \one }\\
 &= \frac 1 n \sum_{i=1}^n \sum_{j=1}^n \ol{w}_{ij} t_{ij} 
  = \frac 1 n \sum_{j=1}^n \ip{W^* T e_j,e_j }  
  =  \frac 1n \Tr(|T|)
\end{align*}
Thus $\| T \|_m = \frac1n \Tr(|T|)$.
\qed

\noindent\textbf{Proof of Theorem~\ref{schur_group}.}
We have already observed that elements of $\T$ are constant on the diagonal.
Thus $\|T\|_m = \frac1n \Tr(|T|) = \ip{|T|e_{x_0},e_{x_0}}$.
For the rest, observe that $W$ belongs to $\ca(T)$.
Hence so does $\ol{W}$ because the basis $T_i$ of $\T$ has real entries.
\qed

We will provide an interesting example in the next section.
For now we provide a couple of more accessible ones.

\begin{eg}\label{off-diag}
Consider the action of the symmetric group $\fS_n$ acting on 
a set $X$ with $n$ elements in the canonical way.
Then the orbits in $X^2$ are just the diagonal $X_0$
and its complement $X_1$.
So $S_{X_1}$ is the projection onto the off-diagonal part of the matrix.

Observe that $X_1 = \one \one^* - I$, where $\one$ is the vector 
of $n$ ones.  
Since $\one\one^* = nP$, where $P$ is the 
projection onto $\bC \one$, $X_1 = (n-1)P - P^\perp$.
Therefore we obtain a formula due to Bhatia, Choi and Davis \cite{BCD}
\begin{align*}
 \| X_1\|_m &= \frac1n \Tr(|X_1|) = \frac1n \Tr\big( (n-1)P + P^\perp \big)\\
 &= \frac1n(n-1+n-1) = 2 - \frac2n .
\end{align*}
\end{eg}

\begin{eg}\label{cyclic}
Consider the cyclic group $C_n$ acting on an $n$-element set, $n\ge3$.
Let $U$ be the unitary operator given by $Ue_k = e_{k+1}$ for 
$1 \le k \le n$, working modulo $n$.
The powers of $U$ yields a basis for the commutant of the group action.

Consider $T = U + I$.
The spectrum of $U$ is just $\{ \omega^k: 0 \le k \le n-1 \}$ where
$\omega = e^{2\pi i/n}$.
Thus the spectrum of $|T|$ consists of the points
\[ |1 + \omega^k| = 2 |\cos(\tfrac{k\pi}n)| \qfor 0 \le k \le n-1 . \]
Hence
\[
 \|T\|_m = \frac1n \Tr(|T|)
 = \frac2n \sum_{k=0}^{n-1} |\cos(\tfrac{k\pi}n)| 
 = \begin{cases}
   \frac{2 \cos(\frac\pi {2n})}{ n \sin(\frac\pi{2n})}
     &n \text{ even}\\
   \frac{2}{ n \sin(\frac\pi {2n})}  &n \text{ odd}
  \end{cases}
\]
Thus the limit as $n$ tends to infinity is $\dfrac4\pi$.
The multiplier norms for the odd cycles decrease to $\frac4\pi$,
while the even cycles increase to the same limit.
\end{eg}

\begin{eg}\label{E:mathias}
Mathias \cite{Ma} considers polynomials in the circulant matrices
$C_z$ given by $C_z e_k = e_{k+1}$ for  $1 \le k < n$
and $C_z e_n = ze_1$, where $|z|=1$.  
This falls into our rubric because there is a 
\textit{diagonal unitary} $D$ so
that $DC_zD^* = wU$ where $U$ is the cycle in the previous example
and $w$ is any $n$th root of $z$.
It is easy to see that conjugation by a diagonal unitary
has no effect on the Schur norm.
Thus any polynomial in $C_z$ is unitarily equivalent to an element
of $\ca(U)$ via the diagonal $D$.
Hence the Schur norm equals the normalized trace of the absolute value.

The most interesting example of this was obtained with $z=-1$ and
$S_n = \sum_{k=0}^{n-1}C_{-1}^k$ which is the matrix with entries 
$\sgn(i-j)$.  So the Schur multiplier defined by $S_n$ is a finite 
Hilbert transform.
Mathias shows that 
\[
 \|S_n\|_m =
 \frac2n \sum_{j=1}^{\lfloor n/2 \rfloor} \cot \frac{(2j-1)\pi}{2n} .
\]
{}From this, he obtains sharper estimates on the norm of triangular truncation
than are obtained in \cite{ACN}.
\end{eg}

\section{Kneser and Johnson Graph Patterns}\label{S:graphs}

In this section, we consider an interesting family of symmetric
patterns which arise commonly in graph theory and combinatorial codes.
The Johnson graphs $J(v,n,i)$ have $\binom v n$ vertices indexed by $n$ element
subsets of a $v$ element set, and edges between $A$ and $B$ if $|A\cap B|=i$.
Thus $0 \le i \le n$.
We consider only $1 \le n \le v/2$ since, if $n>v/2$, one obtains the
same graphs by considering the complementary sets of cardinality $v-n$.
We will explicitly carry out the calculation for the Kneser graphs
$K(v,n) = J(v,n,0)$, and in particular, for $K(2n+1,n)$.
For more on Johnson and Kneser graphs, see \cite{GR}.

We obtained certain Kneser graphs from Toeplitz patterns.
Take a finite subset $S = \{ s_1 < s_2 < \dots < s_{2n+1} \}$
and consider the Toeplitz pattern $\P$ with diagonals in $S$,
namely $\P = \{(i,j) : j-i \in S \}$.
Consider $R$ to be the set of all sums of $n$ elements from $S$
and $C$ to be the set of all sums of $n+1$ elements from $S$.
Index $R$ by the corresponding subset $A$ of $\{1,2,\dots,2n+1\}$
of cardinality $n$; and
likewise index each element of $C$ by a subset $B$ of cardinality $n+1$.
Then for each entry $A$ in $R$, there are exactly $n+1$ elements of $C$
which contain it.  The difference of the sums is an element of $S$.
It is convenient to re-index $C$ by sets of cardinality $n$, replacing $B$ by
its complement $\{1,2,\dots,2n+1\} \bsl B$.
Then the pattern can be seen to be the Kneser graph $K(2n+1,n)$
with $\binom{2n+1}{n}$ vertices indexed by $n$ element subsets of a $2n+1$
element set, with an edge between vertices $A$ and $B$ if $A\cap B = \mt$.
In general, unfortunately, $\P \cap (R \times C)$ will contain more than just these entries,
because two subsets of $S$ of size $n+1$ can have the same sum.

The adjacency matrix of a graph $\G$ is a $v \times v$ matrix with a 1 in each
entry $(i,j)$ corresponding to an edge from vertex $i$ to vertex $j$, and 0's
elsewhere.  
This is a symmetric matrix and its spectral theory is available in the graph 
theory literature; see, for example, \cite{GR}. 
We prove the simple facts we need. 

Fix $(v,n)$ with $n \le v$ and let $X$ denote the set of $n$ element subsets of 
$\{1,\ldots,v\}$.
Define a Hilbert space $\H=\H_X$ as in the previous section
but write the basis as $\{e_A : A \in X \}$.
Observe that there is a natural action $\pi$ of the symmetric group $\fS_v$
on $X$.
The orbits in $X^2$ are 
\[ X_i = \{ (A,B) : A,B \in X,\ |A \cap B| = i \} \qfor 0 \le i \le
n . \] The matrix $T_i$ is just the adjacency matrix of the Johnson
graph $J(v,n,i)$ and, in particular, $T_n = I$.

This action has additional structure that does not
hold for arbitrary transitive actions.

\begin{lem}
The commutant $\T = \spn \{ T_i : 0 \le i \le n \}$ of $\pi(\fS_v)$ is abelian.
Thus $\pi$ decomposes into a direct sum 
of $n+1$ distinct irreducible representations.
\end{lem}

\Prf 
Equality with the span was observed in the last section.
To see that the algebra $\T$ is abelian, observe that
$ T_i T_j = \sum_{k=0}^n a_{ijk} T_k $
where we can find the coefficients $a_{ijk}$ by fixing any two sets
$A,B \subset V$ of size $n$ with $|A \cap B| = k$ and computing
\[ a_{ijk} = |\{ C \subset V: |C|=n,\ |A \cap C|=i,\ |C\cap B| = j \} .\]
This is clearly independent of the order of $i$ and $j$.
As $\T$ is abelian and $n+1$ dimensional, the representation $\pi$
decomposes into a direct sum of $n+1$
distinct irreducible representations.
\qed

\begin{cor}
$\|T_i\| = \binom n i \binom {v-n}{n-i}$ and this is an eigenvalue of
multiplicity one.
The spectrum of $T_i$ contains at most $n+1$ points.  
\end{cor}

\Prf Observe that if $|A|=n$, then the number of subsets $B \in X$
with $|A\cap B | = i$ is $\binom n i \binom {v-n}{n-i}$.
Thus $T_i$ has this many $1$'s in each row.
Hence
\[ T_i \one = \binom n i \binom {v-n}{n-i} \one . \]
Clearly $T_i$ has nonnegative entries and is indecomposable (except for
$i=n$, the identity matrix). 
So by the Perron--Frobenius Theorem, $\binom n i \binom {v-n}{n-i}$ is the
spectral radius and $\one$ is the unique eigenvector;  
and there are no other eigenvalues on the circle of this radius.  
Since $T=T^*$, the norm equals spectral radius.
As $\T$ is $n+1$ dimensional, the spectrum can have at most $n+1$ points.
\qed

We need to identify the invariant subspaces of $\fS_v$ 
as they are the eigenspaces of $T_i$.
The space $V_0 = \bC \one$ yields the trivial representation.
Define vectors associated to sets $C \subseteq \{1,\ldots,v\}$ 
of cardinality at most $n$, including the empty set, by
\[  v_C := \sum_{|A|=n ,\ A \cap C = \mt} e_A .\]
Then define subspaces 
$ V_i = \spn \{ v_C : |C| = i \}$ for $0 \le i \le n .$
It is obvious that each $V_i$ is invariant for $\fS_v$.
Given $C$ with $|C|=i$, we have
\[
 \sum_{C\subset D, |D|=i+1} v_D = (v-n-i) v_C,
\]
as the coefficient of $e_A$ counts the number of choices for the
$(i+1)$st element of $D$ disjoint from an $A$ already disjoint
from $C$.  Therefore 
\[
 \bC \one = V_0 \subset V_1 \subset V_2 \subset \dots \subset V_n .
\]
So the $n+1$ subspaces $W_i = V_i \ominus V_{i-1}$ are invariant for $\fS_v$.

Let $E_i$ denote the idempotent in $\T$ projecting onto $W_i$.
Observe that $\T = \spn\{E_i : 0 \le i \le n\}$.
We need to know the dimension of these subspaces.

\begin{lem}\label{dimeig}
The vectors $\{ v_C : |C|=i \}$ are linearly independent.
Hence $\Dim W_i = \binom v{i} - \binom v{i-1}$.
\end{lem}

\Prf Suppose that $v_{C_0} + \sum_{|C|=i,\ C\ne C_0} \gamma_C v_C = 0$.
By averaging over the subgroup of $\fS_v$ which fixes $C_0$,
namely $\fS_i \times \fS_{v-i}$, we may assume that the coefficients
are invariant under this action.  
Hence $\gamma_C = \alpha_j$ where $j = |C \cap C_0|$.
So with $w_j := \sum_{|C|=i,\ |C \cap C_0| = j} v_C$, we have
$\sum_{j=0}^i \alpha_j w_j = 0$ where $\alpha_i=1$.
We also define vectors $x_k = \sum_{|A \cap C_0| = k} e_A$,
which are clearly linearly independent for $0 \le k \le i$.
Compute for $0 \le j \le i$ (here $A$ implicitly has $|A|=n$)
\[
  w_j = \sum_{\substack{|C|=i\\ |C \cap C_0| = j}} \sum_{A \cap C = \mt} e_A 
  = \sum_{k=0}^{i-j} b_{jk} x_k
\]
where the coefficients are obtained by counting, for a fixed set $A$ with
$|C_0 \cap A| = k$ and $k \le i-j$:
\[
 b_{jk} = |\{C: |C|\!=\!i,\, |C \cap C_0| \!=\! j,\, A \cap C \!=\! \mt\}|
 = \binom {i\!-\!k} j \binom{v\!+\!k\!-\!n\!-\!i}{i\!-\!j} .
\]
It is evident by induction that 
\[
 \spn\{w_j : i-k \le j \le i \} = \spn\{ x_j : 0 \le j \le k \} .
\]
So $\{ v_C : |C|=i \}$ are linearly independent.
\qed

We write $T_i = \sum_{j=0}^n \lambda_{ij} E_j$ be the spectral
decomposition of each $T_i$. 
The discussion above shows that if $|C|=j$, then $v_C$ is 
contained in $V_j$ but not $V_{j-1}$.  
Thus $\lambda_{ij}$ is the unique scalar so that 
$(T_i-\lambda_{ij} I) v_C \in V_{j-1}$. 
This idea can be used to compute the eigenvalues, but the computations
are nontrivial.  
We refer to \cite[Theorem~9.4.3]{GR} for the Kneser graph 
$K(2n+1,n)$ which is the only one we work out in detail.

\begin{lem}\label{eigenlemma}
The adjacency matrix for the Kneser graph $K(2n+1,n)$
has eigenvalues are $(-1)^i(n+1-i)$ with
eigenspaces $W_i$ for $0 \le i \le n$.
\end{lem}

\begin{thm}\label{schur_kneser}
If $T$ is the adjacency matrix of $K(2n+1,n)$, then
\[
 \|T\|_m = \|S_T|_{\T}\| = \dfrac{2^{2n}}{\binom{2n+1}n} = 
 \frac{(4)(6)\dots(2n+2)}{(3)(5)\dots(2n+1)} > \frac12 \log (2n+3) .
\]
\end{thm}

\Prf By Theorem~\ref{schur_group} and Lemma~\ref{eigenlemma},
{\allowdisplaybreaks
\begin{align*}
 \|T\|_m &= \| \Delta(|T|) \| = 
 \binom{2n+1}{n}^{-1} \sum_{i=0}^n (n+1-i) \Tr(E_i) \\
  &= \binom{2n+1}{n}^{-1} \sum_{i=0}^n (n+1-i) 
   \left( \binom{2n+1}{i} - \binom{2n+1}{i-1} \right) \\
  &= \binom{2n+1}{n}^{-1} \sum_{i=0}^n \binom{2n+1}{i} \\
  &= \binom{2n+1}{n}^{-1} \frac12 \sum_{i=0}^{2n+1} \binom{2n+1}{i}
\\&= \binom{2n+1}{n}^{-1} 2^{2n}
   = \dfrac{2^{2n}n! (n+1)!}{(2n+1)!} \\&= 
  \dfrac{2\cdot 4 \cdots(2n)\,\,\, 2\cdot 4 \cdots(2n)\cdot(2n+2)}
  {2\cdot 4 \cdots(2n)\,\,1\cdot3\cdots
  (2n\!-\!1)(2n\!+\!1)} \\
  &= \dfrac{2\cdot 4 \cdots(2n)\cdot(2n+2)}
  {1\cdot3\cdots (2n\!-\!1)(2n\!+\!1)} \\
  &= \prod_{i=0}^n \left( 1 + \frac1{2i+1} \right) > \frac12 \log (2n+3) .
\end{align*}
} 
\upqed


\end{document}